\documentclass[11pt]{article}

\usepackage{amssymb}
\usepackage{amsfonts}
\usepackage{amsmath}
\usepackage{amsthm}
\usepackage[margin=1.00in]{geometry}
\usepackage{enumerate}
\usepackage{amstext}
\usepackage{layout}

\numberwithin{equation}{section}

\newtheorem{theorem}{Theorem}[section]

\newtheorem{lemma}{Lemma}[section]

\newtheorem{remark}{Remark}[section]

\begin{document}

\begin{center}
{\large \bf An Extension of Feller's Strong Law of Large Numbers}

\vskip 0.5cm

{\bf Deli Li$^{1}$, 
Han-Ying Liang$^{2,}$\footnote{Corresponding author: Han-Ying Liang (Telephone: 86-21-65983242, FAX: 86-21-65983242)}, 
and Andrew Rosalsky$^{3}$}

\vskip 0.3cm

$^1$Department of Mathematical Sciences, Lakehead
University,\\
Thunder Bay, Ontario, Canada\\
$^2$Department of Mathematics, Tongji University,\\
Shanghai, China\\
$^3$Department of Statistics, University of Florida,\\
Gainesville, Florida, USA

\end{center}

\vskip 0.3cm

\noindent {\bf Abstract}~~This paper presents a general result that allows for establishing a link between
the Kolmogorov-Marcinkiewicz-Zygmund strong law of large numbers and Feller's strong law of large numbers 
in a Banach space setting. Let $\{X, X_{n}; n \geq 1\}$ be a sequence of independent and identically distributed 
Banach space valued random variables and set $S_{n} = \sum_{i=1}^{n}X_{i},~n \geq 1$. Let 
$\{a_{n}; n \geq 1\}$ and $\{b_{n}; n \geq 1\}$ be increasing sequences of positive real numbers such that 
$\lim_{n \rightarrow \infty} a_{n} = \infty$ and $\left\{b_{n}/a_{n};~ n \geq 1 \right\}$ is a nondecreasing 
sequence. We show that 
\[
\frac{S_{n}- n \mathbb{E}\left(XI\{\|X\| \leq b_{n} \} \right)}{b_{n}} \rightarrow 0~~\mbox{almost surely}
\] 
for every Banach space valued random variable $X$ with $\sum_{n=1}^{\infty} \mathbb{P}(\|X\| > b_{n}) < \infty$ 
if $S_{n}/a_{n} \rightarrow 0$ almost surely for every symmetric Banach space valued random variable $X$ with 
$\sum_{n=1}^{\infty} \mathbb{P}(\|X\| > a_{n}) < \infty$. To establish this result, we invoke two tools
(obtained recently by Li, Liang, and Rosalsky): a symmetrization procedure for the strong law of large numbers 
and a probability inequality for sums of independent Banach space valued random variables.

~\\


\vskip 0.3cm


\vskip 0.3cm

\noindent {\bf Running Head}: On Feller's Strong Law of Large Numbers

\section{Introduction and the main result}

We begin with stating Feller's (1946) strong law of large numbers (SLLN) as follows.

\vskip 0.3cm

\noindent {\bf Theorem A}.~{\rm (Feller's SLLN. Theorems 1 and 2 of Feller (1946))}.
{\it Let $\{X, X_{n}; n \geq 1 \}$ be a sequence of independent and identically distributed 
(i.i.d.) real-valued random variables, and let $S_{n} = \sum_{i=1}^{n}$, $n \geq 1$. 
Let $\{b_{n}; n \geq 1 \}$ be an increasing sequence of positive real numbers.
Suppose that one of the following two sets of conditions holds:

\noindent {\bf (i)}~~For some $0 < \delta < 1$, $\mathbb{E}|X|^{1 + \delta} = \infty$, $\mathbb{E}X = 0$,
and there exists an $\epsilon$ with $0 \leq \epsilon < 1$
such that
\[
b_{n}n^{-1/(1 + \epsilon)} \uparrow ~~\mbox{and}~~b_{n}/n \downarrow,
\]
{\bf (ii)}~~$\mathbb{E}|X| = \infty$ and
\[
b_{n}/n \uparrow.
\]
Then we have 
\[
\lim_{n \rightarrow \infty} \frac{S_{n}}{b_{n}} 
= 0~~\mbox{almost surely (a.s.) or}~~\limsup_{n \rightarrow \infty} \frac{\left|S_{n}\right|}{b_{n}} 
= \infty~~\mbox{a.s.}
\]
according as
\[
\sum_{n=1}^{\infty}\mathbb{P}(|X| > b_{n}) < \infty ~~\mbox{or}~~ = \infty.
\]
}

Feller's SLLN is a remarkable limit theorem concerning sums of i.i.d. random variables. 
From the internet, one can find that Feller's SLLN has received more than 150 citations
where many of them have been received within the most recent 5 years.

This paper presents a general result in a Banach space setting that allows for establishing a link between
the Kolmogorov-Marcinkiewicz-Zygmund SLLN and Feller's SLLN.

For stating our main result, we introduce some notation as follows. Let $(\Omega, \mathcal{F}, \mathbb{P})$ 
be a probability space and let $(\mathbf{B}, \| \cdot \| )$ be a real separable Banach space equipped with 
its Borel $\sigma$-algebra $\mathcal{B}$ ($=$ the $\sigma$-algebra generated by the class of open subsets 
of $\mathbf{B}$ determined by $\|\cdot\|$). A {\bf B}-valued random variable $X$ is defined as 
a measurable function from $(\Omega, \mathcal{F})$ into $(\mathbf{B}, \mathcal{B})$. 
Let $ \{X, X_{n};~n \geq 1 \}$  be a sequence of i.i.d. {\bf B}-valued random variables and put 
$S_{n} = \sum_{i=1}^{n} X_{i}$, $n \geq 1$. Let $\{R, R_{n}; ~n \geq 1\}$ be a {\it Rademacher sequence}; 
that is, $\{R_{n};~n \geq 1\}$ is a sequence of i.i.d. random variables 
with $\mathbb{P}\left(R = 1\right) = \mathbb{P}\left(R = -1\right) = 1/2$. Let
$\mathbf{B}^{\infty} =
\mathbf{B}\times\mathbf{B}\times\mathbf{B}\times \cdots$ and define
\[
\mathcal{C}(\mathbf{B}) = \left\{(v_{1}, v_{2}, ...) \in
\mathbf{B}^{\infty}: ~\sum_{n=1}^{\infty} R_{n}v_{n}
~~\mbox{converges in probability} \right\}.
\]
Let $1 \leq p \leq 2$. Then $\mathbf{B}$ is said to be of {\it
Rademacher type $p$} if there exists a constant $0 < C < \infty$
such that
\[
\mathbb{E}\left\|\sum_{n=1}^{\infty}R_{n}v_{n} \right\|^{p} \leq C
\sum_{n=1}^{\infty}\|v_{n}\|^{p}~~\mbox{for all}~(v_{1}, v_{2}, ...)
\in \mathcal{C}(\mathbf{B}).
\]
The following remarkable theorem, which is due to de Acosta (1981), provides 
a characterization of Rademacher type $p$ Banach spaces.

\vskip 0.3cm

\noindent {\bf Theorem B}.~{\rm (de Acosta (1981))}. 
{\it Let $1 \leq p < 2$. Then the following two
statements are equivalent:
\begin{align*}
& {\bf (i)} \quad \mbox{The Banach space $\mathbf{B}$ is of
Rademacher type $p$.}\\
& {\bf (ii)} \quad \mbox{For every sequence $\{X, X_{n}; ~n \geq 1 \}$
of i.i.d. {\bf B}-valued variables},
\end{align*}
\[
\lim_{n \rightarrow \infty} \frac{S_{n}}{n^{1/p}} = 0~~\mbox{a.s. if
and only if}~~\mathbb{E}\|X\|^{p} < \infty~~\mbox{and}~~\mathbb{E}X
= 0.
\]
}

The main result of this paper is the following theorem. 

\vskip 0.3cm

\begin{theorem}
Let $(\mathbf{B}, \|\cdot\|)$ be a real separable Banach space.
Let $\{a_{n}; n \geq 1\}$ and $\{b_{n}; n \geq 1\}$ be increasing
sequences of positive real numbers such that 
\begin{equation}
\lim_{n \rightarrow \infty} a_{n} = \infty~ \mbox{and}~ b_{n}/a_{n} \uparrow.
\end{equation}
Suppose that for every symmetric sequence
$\{X, X_{n}; ~n \geq 1 \}$ of i.i.d. {\bf B}-valued random variables,
\begin{equation}
\lim_{n \rightarrow \infty} \frac{S_{n}}{a_{n}} = 0 
~~\mbox{a.s. whenever}~~\sum_{n=1}^{\infty}\mathbb{P}(\|X\| > a_{n}) < \infty.
\end{equation}
Then, for every sequence $\{X, X_{n}; ~n \geq 1 \}$ of i.i.d. {\bf B}-valued random variables, we have that
\begin{equation}
\lim_{n \rightarrow \infty} \frac{S_{n}- \gamma_{n}}{b_{n}} 
= 0~~\mbox{a.s. or}~~\limsup_{n \rightarrow \infty} \frac{\left\|S_{n}- \gamma_{n} \right\|}{b_{n}} 
= \infty~~\mbox{a.s.}
\end{equation}
according as
\begin{equation}
\sum_{n=1}^{\infty}\mathbb{P}(\|X\| > b_{n}) < \infty ~~\mbox{or}~~ = \infty.
\end{equation}
Here and below $\gamma_{n} = n \mathbb{E}\left(XI\{\|X\| \leq b_{n} \} \right)$, $n \geq 1$.
\end{theorem}

\vskip 0.3cm

\begin{remark}
We now can see how Feller's SLLN can be easily derived from the Kolmogorov-Marcinkiewicz-Zygmund SLLN and Theorem 1.1 above.
For given $\epsilon$ with $0 \leq \epsilon < 1$, write $a_{n} = n^{1/(1 + \epsilon)},~n \geq 1$. The celebrated 
Kolmogorov-Marcinkiewicz-Zygmund SLLN (see Kolmogoroff (1930) for $\epsilon = 0$ and Marcinkiewicz and Zygmund 
(1937) for $0 < \epsilon < 1$) asserts that, for every symmetric sequence $\{X, X_{n}; ~n \geq 1 \}$ of i.i.d. real-valued
(i.e., ${\bf B} = \mathbb{R}$) random variables
\[
\lim_{n \rightarrow \infty} \frac{S_{n}}{a_{n}} = \lim_{n \rightarrow \infty} \frac{S_{n}}{n^{1/(1 + \epsilon)}} = 0 
~~\mbox{a.s. if and only if}~~\mathbb{E}|X|^{1+ \epsilon} < \infty
~~(\mbox{i.e.},~\sum_{n=1}^{\infty}\mathbb{P}(|X| > a_{n}) < \infty).
\]
Then by Theorem 1.1, for every sequence $\{X, X_{n}; ~n \geq 1 \}$ of i.i.d. real-valued random variables
and every increasing sequence $\{b_{n}; n \geq 1 \}$ of positive real numbers with 
\[
b_{n}/a_{n} = b_{n} n^{-1/(1+ \epsilon)} \uparrow,
\]
we have
\[
\lim_{n \rightarrow \infty} \frac{S_{n}- \gamma_{n}}{b_{n}} 
= 0~~\mbox{a.s. or}~~\limsup_{n \rightarrow \infty} \frac{\left |S_{n}- \gamma_{n} \right|}{b_{n}} 
= \infty~~\mbox{a.s.}
\]
according as
\[
\sum_{n=1}^{\infty}\mathbb{P}(|X| > b_{n}) < \infty ~~\mbox{or}~~ = \infty,
\]
where $\gamma_{n} = n \mathbb{E}\left(XI\{|X| \leq b_{n} \} \right)$, $n \geq 1$. That is, Feller's SLLN follows
from the Kolmogorov-Marcinkiewicz-Zygmund SLLN and Theorem 1.1 above.
\end{remark}

\vskip 0.3cm

\begin{remark}
Under the assumptions of Theorem 1.1, it follows from the conclusion of Theorem 1.1 that, for every 
sequence $\{X, X_{n}; ~n \geq 1 \}$ of i.i.d. {\bf B}-valued random variables,
\[
\lim_{n \rightarrow \infty} \frac{S_{n}- \gamma_{n}}{b_{n}} 
= 0~~\mbox{a.s. if and only if}~~\sum_{n=1}^{\infty}\mathbb{P}(\|X\| > b_{n}) < \infty,
\]
\[
\limsup_{n \rightarrow \infty} \frac{\|S_{n}- \gamma_{n}\|}{b_{n}} 
= \infty~~\mbox{a.s. if and only if}~~\sum_{n=1}^{\infty}\mathbb{P}(\|X\| > b_{n}) = \infty.
\]
Hence under the assumptions of Theorem 1.1, there does not exist a sequence $\{X, X_{n}; ~n \geq 1 \}$ of i.i.d. 
{\bf B}-valued random variables such that
\[
0 < \limsup_{n \rightarrow \infty}
\frac{\left\|S_{n} - \gamma_{n} \right\|}{b_{n}} < \infty~~\mbox{a.s.}
\]
\end{remark}

\vskip 0.3cm

Also, combining Theorem 1.1 and Theorem B above, we immediately obtain the following two remarks.

\vskip 0.3cm

\begin{remark}
Let $1 \leq p < 2$ and let $\{a_{n}; n \geq 1\}$ be an increasing sequences of positive 
real numbers such that 
\[
\lim_{n \rightarrow \infty} a_{n} = \infty~ \mbox{and}~ n^{1/p}/a_{n} \uparrow.
\]
Let $(\mathbf{B}, \|\cdot\|)$ be a real separable Banach space such that (1.2) holds for every 
symmetric sequence $\{X, X_{n}; n \geq 1 \}$ of i.i.d. {\bf B}-valued random variables. Then 
the Banach space ${\bf B}$ is of Rademacher type $p$.  
\end{remark}

\vskip 0.3cm

\begin{remark}
Let $1 \leq p < 2$ and let $\{b_{n}; n \geq 1\}$ be a sequence of positive 
real numbers such that
\[
b_{n}/n^{1/p} \uparrow.
\] 
If $\mathbf{B}$ is of Rademacher type $p$, then for every sequence 
$\{X, X_{n}; ~n \geq 1 \}$ of i.i.d. {\bf B}-valued random variables, 
(1.3) and (1.4) are equivalent.
\end{remark}

\vskip 0.3cm

\begin{remark}
Remark 1.4 should be compared with Theorem 4 of Adler, Rosalsky, and Taylor (1989) and with 
(in the unweighted case) the key lemma (Lemma 6) of that article. In Lemma 6 (wherein 
$1 \leq p \leq 2$) and Theorem 4 (wherein $1 < p \leq 2$) of Adler, Rosalsky, and Taylor (1989),
for a sequence $\left \{X, X_{n}; ~n \geq 1 \right \}$ of i.i.d. random variables in 
a real separable Rademacher type $p$ Banach space and a sequence of positive constants 
$b_{n} \uparrow \infty$, conditions are provided under which
\[
\sum_{n=1}^{\infty} \mathbb{P} \left(\|X\| > b_{n} \right) < \infty
\]
ensures that
\[
\lim_{n \rightarrow \infty} \frac{\sum_{i=1}^{n} \left(X_{i} - \mathbb{E} \left(XI\{\|X\| \leq b_{i}\} \right) \right)}{b_{n}}
= 0~~\mbox{a.s.}
\]
and 
\[
\lim_{n \rightarrow \infty} \frac{\sum_{i=1}^{n} X_{i}}{b_{n}} = 0~~\mbox{a.s.}
\]
\end{remark}

\vskip 0.3cm

The proof of Theorem 1.1 will be provided in Section 2. To establish Theorem 1.1,
we invoke two tools (obtained recently by Li, Liang, and Rosalsky (2017a, b)): a 
symmetrization procedure for the SLLN and a probability inequality which is a
comparison theorem for sums of independent ${\bf B}$-valued random variables.

\vskip 0.3cm

We close this section by remarking that a version of Feller's SLLN was obtained by
Martikainen and Petrov (1980) for a sequence of identically distributed real-valued
random variables $\left \{X, X_{n};~n \geq 1 \right \}$ without any independence
conditions being imposed on the summands; the result holds irrespective of the joint 
distributions of the summands. Specifically, in Theorem 2 of Martikainen and Petrov 
(1980) it is shown that
\[
\lim_{n \rightarrow \infty} \frac{\sum_{i=1}^{n} X_{i}}{b_{n}} = 0~~\mbox{a.s.}
\]
if 
\[
0 < b_{n} \uparrow, ~~b_{n} \sum_{i=n}^{\infty} \frac{1}{b_{i}} = {\it O}(n),~~\mbox{and}~~
\sum_{n=1}^{\infty} \mathbb{P} \left(|X| > b_{n} \right) < \infty.
\]

\section{Proof of Theorem 1.1}

Throughout this section, $\{a_{n}; n \geq 1\}$ and $\{b_{n}; n \geq 1\}$ are increasing
sequences of positive real numbers satisfying (1.1). Write
\[
I(1) = \left \{i;~ b_{i} \leq 2 \right \} ~~\mbox{and}~~
I(m) = \left \{i;~2^{m-1} < b_{i} \leq 2^{m} \right\},~~m \geq 2.
\]
It follows from (1.1) that $0 < b_{n} \uparrow \infty$ and
\[
\{1, 2, 3, ..., n, ... \} = \bigcup_{m=1}^{\infty} I(m).
\]
Note that the $I(m), ~m \geq 1$ are mutually exclusive sets. 
Thus there exist positive integers $k_{n}, m_{n}, ~n \geq 1$ such that
\[
k_{1} < k_{2} < ... < k_{n} < ..., ~~m_{1} < m_{2} < ... < m_{n} < ...,
\]
\[
\{1, 2, 3, ..., n, ... \} = \bigcup_{n=1}^{\infty} I\left(m_{n} \right), ~~\mbox{and}~~I\left(m_{1} \right) 
= \left\{1, ..., k_{1} \right\}, ~I\left(m_{n} \right) = \left\{k_{n-1} + 1, ..., k_{n} \right\}, ~n \geq 2.
\]

To prove Theorem 1.1, we use the following four preliminary lemmas. 

\vskip 0.3cm

\begin{lemma}
{\rm (Lemma 3.1 of Li, Liang, and Rosalsky (2017 b))}~
There exist two continuous and increasing functions $\varphi(\cdot)$ and $\psi(\cdot)$
defined on $[0, \infty)$ such that 
\begin{equation}
\lim_{t \rightarrow \infty} \varphi(t) = \infty 
~~\mbox{and}~~\frac{\psi(\cdot)}{\varphi(\cdot)} ~\mbox{is a nondecreasing function on}~[0, \infty),
\end{equation}
\begin{equation}
\varphi(0) = \psi(0) = 0, ~\varphi(n) = a_{n}, ~\psi(n) = b_{n}, ~ n \geq 1.
\end{equation}
\end{lemma}

\vskip 0.3cm

\begin{lemma}
Let $\{V_{n};~n \geq 1 \}$ be a sequence of independent and symmetric {\bf B}-valued random variables.
Set $k_{0} = 0$. Then the following two statements hold.

{\bf (i)} ~~If 
\begin{equation}
\lim_{n \rightarrow \infty} \frac{\sum_{i=1}^{n}V_{i}}{a_{n}} = 0~~\mbox{a.s.,} 
\end{equation}
then
\begin{equation}
\sum_{n=1}^{\infty} 
\mathbb{P}\left(\left\|\sum_{i = k_{n-1} + 1}^{k_{n}} V_{i} \right\| > \epsilon a_{k_{n}} \right)
< \infty ~~\forall~ \epsilon > 0.
\end{equation}

{\bf (ii)}
\begin{equation}
\lim_{n \rightarrow \infty} \frac{\sum_{i=1}^{n}V_{i}}{b_{n}} = 0~~\mbox{a.s.}
\end{equation}
if and only if
\begin{equation}
\sum_{n=1}^{\infty} 
\mathbb{P}\left(\left\|\sum_{i=k_{n-1}+1}^{k_{n}} V_{i} \right\| > \epsilon b_{k_{n}} \right)
< \infty ~~\forall~ \epsilon > 0.
\end{equation}
\end{lemma}

\noindent {\it Proof}~~We first prove Part (i). Clearly, (2.3) implies that
\begin{equation}
\frac{\left\|\sum_{i = k_{n-1} + 1}^{k_{n}} V_{i} \right\|}{a_{k_{n}}} 
\leq \frac{\left\|\sum_{i = 1}^{k_{n}} V_{i} \right\|}{a_{k_{n}}} + \left(\frac{a_{k_{n-1}}}{a_{k_{n}}} \right)
\frac{\left\|\sum_{i = 1}^{k_{n-1}} V_{i} \right\|}{a_{k_{n-1}}} \rightarrow 0 ~~\mbox{a.s. as}~~n \rightarrow \infty.
\end{equation}
Since the $\sum_{i=k_{n-1} + 1}^{k_{n}} V_{i}, ~n \geq 1$ are independent {\bf B}-valued random variables, (2.4) follows 
from (2.7) and the Borel-Cantelli lemma.

We now establish Part (ii). From the proof of Part (i), we only need to show that (2.5) follows from (2.6).
Since $\{V_{n};~n \geq 1 \}$ is a sequence of independent and symmetric {\bf B}-valued random variables,
by the remarkable L\'{e}vy inequality in a Banach space setting (see, e.g., see Proposition 2.3 of Ledoux and Talagrand 
(1991)), we have that for every $n \geq 1$,
\[
\mathbb{P} \left( \max_{k_{n-1} < k \leq k_{n}}
\left \|\sum_{i= k_{n-1} + 1}^{k} V_{i} \right \| > \epsilon b_{k_{n}}  \right) 
\leq 2 \mathbb{P}\left(\left\|\sum_{i = k_{n-1} + 1}^{k_{n}}V_{i} \right\| > \epsilon b_{k_{n}} \right)
~~\forall~\epsilon > 0.
\]
Thus it follows from (2.6) that
\[
\sum_{n=1}^{\infty} \mathbb{P} \left( \max_{k_{n-1} < k \leq k_{n}}
\left \|\sum_{i= k_{n-1} + 1}^{k} V_{i} \right \| > \epsilon b_{k_{n}}  \right) < \infty 
~~\forall~\epsilon > 0
\]
which ensures that
\begin{equation}
A_{n} \stackrel{\Delta}{=}
\frac{\max_{k_{n-1} < k \leq k_{n}} \left \|\sum_{i= k_{n-1} + 1}^{k} V_{i} \right \|}{b_{k_{n}}}
\rightarrow 0 ~~\mbox{a.s.}
\end{equation}
Now by the Toeplitz lemma, we conclude from (2.8) that
\[
\begin{array}{lll}
\mbox{$\displaystyle 
\max_{k_{n-1} < k \leq k_{n}} \frac{\left \|\sum_{i= 1}^{k} V_{i} \right \|}{b_{k}} $}
&\leq& 
\mbox{$\displaystyle 
2 \max_{k_{n-1} < k \leq k_{n}} \frac{\left \|\sum_{i= 1}^{k} V_{i} \right \|}{b_{k_{n}}} $}\\
&&\\
&\leq&
\mbox{$\displaystyle 
2 \left(\sum_{j=1}^{n-1} \frac{\left \|\sum_{i= k_{j-1} + 1}^{k_{j}} V_{i} \right \|}{b_{k_{n}}}
+ A_{n} \right)$}\\
&&\\
&\leq&
\mbox{$\displaystyle 
2 \sum_{j=1}^{n} \left(\frac{b_{k_{j}}}{b_{k_{n}}} \right) A_{j} $}\\
&&\\
&\leq& 
\mbox{$\displaystyle 
4 \sum_{j=1}^{n} \left(\frac{2^{m_{j}}}{2^{m_{n}}} \right) A_{j} $} \\
&&\\
&\rightarrow& 
\mbox{$\displaystyle 
0 ~~\mbox{a.s.~as}~ n \rightarrow \infty;$}
\end{array}
\]
i.e., (2.5) holds. ~$\Box$

\vskip 0.3cm

The following probability inequality is due to Li, Liang, and Rosalsky (2017 b)
and is a comparison theorem for sums of independent $\mathbf{B}$-valued random variables.

\vskip 0.2cm

\begin{lemma}
{\rm (Theorem 1.1 (ii) of Li, Liang, and Rosalsky (2017 b))}~
Let $\varphi(\cdot)$ and $\psi(\cdot)$ be two continuous and increasing functions
defined on $[0, \infty)$ with $\varphi(0) = \psi(0) = 0$ and satisfying (2.1).  
If $\{V_{n};~n \geq 1 \}$ is a sequence of independent 
and symmetric {\bf B}-valued random variables, then for every 
$n \geq 1$ and all $t \geq 0$,
\[
\mathbb{P}\left(\left\|\sum_{i=1}^{n} V_{i} \right\| > t b_{n} \right) 
\leq 4 \mathbb{P} \left(\left\|\sum_{i=1}^{n} \varphi\left(\psi^{-1}(\|V_{i}\|)\right) 
\frac{V_{i}}{\|V_{i}\|} \right\| > t a_{n} \right) 
+ \sum_{i=1}^{n}\mathbb{P}\left(\|V_{i}\| > b_{n} \right).
\]
\end{lemma}

\vskip 0.3cm

The following symmetrization procedure for the SLLN for independent ${\bf B}$-valued random variables
is due to Li, Liang, and Rosalsky (2017 a).

\vskip 0.3cm

\begin{lemma}
{\rm (Corollary 1.3 of Li, Liang, and Rosalsky (2017 a))}
Let $\{X, X_{n}; ~n \geq 1\}$ be a sequence of i.i.d. 
{\bf B}-valued random variables. Let $\{X_{n}^{\prime};~n \geq 1 \}$ be an
independent copy of $\{X_{n};~n \geq 1 \}$. Write $S_{n} = \sum_{i=1}^{n} X_{i}$,
$S_{n}^{\prime} = \sum_{i=1}^{n} X_{i}^{\prime}$, $n \geq 1$. Let $\{b_{n}; n \geq 1\}$ 
be an increasing sequence of positive real numbers such that 
$\lim_{n \rightarrow \infty} b_{n} = \infty$. Then 
\[
\lim_{n \rightarrow \infty} 
\frac{S_{n}- n \mathbb{E}\left(XI\{\|X\| \leq b_{n} \} \right)}{b_{n}} = 0 ~~\mbox{a.s.}
\]
if and only if 
\[
\lim_{n \rightarrow \infty} \frac{S_{n} - S_{n}^{\prime}}{b_{n}} = 0~~\mbox{a.s.} 
\]
\end{lemma}

\vskip 0.3cm

With the preliminaries accounted for, Theorem 1.1 may be proved.

\vskip 0.3cm

\noindent {\it Proof of Theorem 1.1}~~To establish this theorem, it suffices to show that, 
for every sequence $\{X, X_{n}; ~n \geq 1 \}$ of i.i.d. {\bf B}-valued random variables, the 
following three statements are equivalent:
\begin{equation}
\lim_{n \rightarrow \infty} \frac{S_{n}- \gamma_{n}}{b_{n}} = 0~~\mbox{a.s.,}
\end{equation}
\begin{equation}
\limsup_{n \rightarrow \infty} \frac{\left\|S_{n}- \gamma_{n} \right\|}{b_{n}} 
< \infty~~\mbox{a.s.,}
\end{equation}
\begin{equation}
\sum_{n=1}^{\infty}\mathbb{P}(\|X\| > b_{n}) < \infty.
\end{equation}

The three statements (2.9)-(2.11) are equivalent if we can show that (2.9) and (2.11) are 
equivalent and (2.9) and (2.10) are equivalent.

For establishing the implication ``(2.9) $\Rightarrow$ (2.11)", let $\{X, X_{n};~n \geq 1 \}$ be a 
sequence of i.i.d. {\bf B}-valued random variables satisfying (2.9). It follows from (2.9) that
\[
\lim_{n \rightarrow \infty} \frac{\sum_{i=1}^{n} \left(X_{i} - X_{i}^{\prime} \right)}{b_{n}} = 0
~~\mbox{a.s.}
\]
which implies that
\begin{equation}
\frac{X_{n} - X_{n}^{\prime}}{b_{n}} \rightarrow 0 ~~\mbox{a.s.}
\end{equation}
By the Borel-Cantelli lemma, (2.12) is equivalent to 
\begin{equation}
\sum_{n=1}^{\infty} \mathbb{P} \left(\|X - X^{\prime}\| > \epsilon b_{n} \right) 
= \sum_{n=1}^{\infty} \mathbb{P} \left(\|X_{n} - X_{n}^{\prime}\| > \epsilon b_{n} \right) 
< \infty ~~\forall ~\epsilon > 0.
\end{equation}
Note that $\left \{\|X^{\prime}\| \leq b_{n}/2, \|X\| > b_{n} \right \} 
\subseteq \left \{ \|X - X^{\prime} \| > b_{n}/2 \right \}$ and
$\displaystyle \lim_{n \rightarrow \infty} \mathbb{P} \left(\|X^{\prime}\| \leq b_{n}/2 \right) = 1$.
We thus have that, for all large $n$,
\[
\mathbb{P} \left(\|X\| > b_{n} \right) \leq 2 \mathbb{P} \left ( \|X - X^{\prime} \| > b_{n}/2 \right ) 
\]
which, together with (2.13), implies (2.11). 

We now prove ``(2.11) $\Rightarrow$ (2.9)". Let $\{X, X_{n}; ~n \geq 1 \}$ be a sequence of i.i.d. 
{\bf B}-valued random variables satisfying (2.11). Since $\{a_{n}; n \geq 1\}$ and $\{b_{n}; n \geq 1\}$ 
are increasing sequences of positive real numbers satisfying (1.1), by Lemma 2.1, there exist
two continuous and increasing functions $\varphi(\cdot)$ and $\psi(\cdot)$
defined on $[0, \infty)$ such that both (2.1) and (2.2) hold. Write
\[
\tilde{X} = \frac{X - X^{\prime}}{2}, ~\tilde{X}_{n} = \frac{X_{n} - X^{\prime}_{n}}{2}, ~n \geq 1
\]
and
\[
Y = \varphi\left(\psi^{-1}(\|\tilde{X}\|)\right) 
\frac{\tilde{X}}{\|\tilde{X}\|}, ~ Y_{n} = \varphi\left(\psi^{-1}(\|\tilde{X}_{n}\|)\right) 
\frac{\tilde{X}_{n}}{\|\tilde{X}_{n}\|}, ~n \geq 1.
\]
Then $\{\tilde{X}, \tilde{X}_{n}; ~n \geq 1 \}$ is
a sequence of i.i.d. symmetric {\bf B}-valued random variables such that
\begin{equation}
\sum_{n=1}^{\infty} \mathbb{P} \left(\|\tilde{X}\| > b_{n} \right) \leq
2 \sum_{n=1}^{\infty} \mathbb{P} \left(\|X\| > b_{n} \right) < \infty
\end{equation}
and $\{Y, Y_{n};~ n\geq 1\}$ is a sequence of of i.i.d. symmetric {\bf B}-valued 
random variables such that
\begin{equation}
\sum_{n=1}^{\infty} \mathbb{P} \left(\|Y\| > a_{n} \right)
= \sum_{n=1}^{\infty} \mathbb{P} \left(\| \tilde{X} \| > b_{n} \right)
< \infty.
\end{equation}
Hence, by Lemma 2.2 (i), we conclude from (2.15) and (1.2) that 
\begin{equation}
\sum_{n=1}^{\infty} \mathbb{P}\left(\left\|\sum_{i = k_{n-1}+1}^{k_{n}}
Y_{i} \right\| > \epsilon a_{k_{n}} \right)
< \infty ~~\forall~ \epsilon > 0.
\end{equation}
By Lemma 2.3, we have that, for every $n \geq 1$,
\[
\begin{array}{lll}
\mbox{$\displaystyle
\mathbb{P}\left(\left\|\sum_{i=k_{n-1}+1}^{k_{n}} \tilde{X}_{i} \right\| > \epsilon b_{k_{n}} \right)$}
& \leq &
\mbox{$\displaystyle  
4 \mathbb{P} \left(\left\|\sum_{i=k_{n-1}+1}^{k_{n}} Y_{i} \right\| > \epsilon a_{n} \right) 
+ \sum_{i=k_{n-1} + 1}^{k_{n}}\mathbb{P}\left(\|\tilde{X}_{i}\| > b_{k_{n}} \right) $}\\
&&\\
&\leq& 
\mbox{$\displaystyle  
4 \mathbb{P} \left(\left\|\sum_{i=k_{n-1}+1}^{k_{n}} Y_{i} \right\| > \epsilon a_{n} \right) 
+ \sum_{i=k_{n-1} + 1}^{k_{n}}\mathbb{P}\left(\|\tilde{X}\| > b_{i} \right)~~\forall~\epsilon > 0.$}
\end{array}
\]
It thus follows from (2.15) and (2.16) that
\begin{equation}
\begin{array}{ll}
& \mbox{$\displaystyle
\sum_{n=1}^{\infty} \mathbb{P}\left(\left\|\sum_{i=k_{n-1}+1}^{k_{n}} \tilde{X}_{i} \right\| > \epsilon b_{k_{n}} \right)$}\\
&\\
& \mbox{$\displaystyle  
\leq 4 \sum_{n=1}^{\infty} \mathbb{P} \left(\left\|\sum_{i=k_{n-1}+1}^{k_{n}} Y_{i} \right\| > \epsilon a_{n} \right) 
+ \sum_{n=1}^{\infty} \sum_{i=k_{n-1} + 1}^{k_{n}}\mathbb{P}\left(\|\tilde{X}_{i}\| > b_{i} \right) $}\\
&\\
& \mbox{$\displaystyle  
= 4 \sum_{n=1}^{\infty} \mathbb{P} \left(\left\|\sum_{i=k_{n-1}+1}^{k_{n}} Y_{i} \right\| > \epsilon a_{n} \right) 
+ \sum_{n=1}^{\infty} \mathbb{P}\left(\|\tilde{X}_{i}\| > b_{n} \right) $}\\
&\\
& \mbox{$\displaystyle  
< \infty~~\forall ~\epsilon > 0.$}
\end{array}
\end{equation}
By Lemma 2.2 (ii), (2.17) is equivalent to
\[
\frac{S_{n} - S_{n}^{\prime}}{2b_{n}} = \frac{\sum_{i=1}^{n} \tilde{X}_{i}}{b_{n}} \rightarrow 0~~\mbox{a.s.}
\]
Hence 
\[
\frac{S_{n} - S_{n}^{\prime}}{b_{n}} \rightarrow 0~~\mbox{a.s.}
\]
By Lemma 2.4, (2.9) follows. 

The implication ``(2.9) $\Rightarrow$ (2.10)" is obvious.

We now establish the implication ``(2.10) $\Rightarrow$ (2.9)". It follows from (2.10) that
\[
\limsup_{n \rightarrow \infty} \frac{\left\|\sum_{i=1}^{n} \left(X_{i} - X_{i}^{\prime} \right) \right\|}{b_{n}} 
< \infty
~~\mbox{a.s.}
\]
which implies that
\begin{equation}
\limsup_{n \rightarrow \infty} \frac{\left\|X_{n} - X_{n}^{\prime} \right\|}{b_{n}} < \infty ~~\mbox{a.s.}
\end{equation}
By the Borel-Cantelli lemma, (2.18) is equivalent to: for some constant $0 < \lambda < \infty$,
\[
\sum_{n=1}^{\infty} \mathbb{P} \left(\|X_{n} - X_{n}^{\prime}\| > \lambda b_{n} \right) 
< \infty; ~~\mbox{i.e.,}~~
\sum_{n=1}^{\infty} \mathbb{P} \left(\left\|\frac{X - X^{\prime}}{\lambda} \right \| > b_{n} \right) 
< \infty.
\]
That is, (2.11) holds with $X$ replaced by symmetric random variable $(X - X^{\prime})/\lambda$. Since (2.9) and (2.11)
are equivalent, we conclude that
\[
\left(\frac{1}{\lambda} \right) \frac{S_{n} - S_{n}^{\prime}}{b_{n}} 
= \frac{\sum_{i=1}^{n} \frac{X_{i} - X_{i}^{\prime}}{\lambda}}{b_{n}} \rightarrow 0 ~~\mbox{a.s.}
\]
Thus
\[
\frac{S_{n} - S_{n}^{\prime}}{b_{n}} \rightarrow 0 ~~\mbox{a.s.}
\]
which, by Lemma 2.4, implies (2.9). The proof of Theorem 1.1 is therefore complete. ~$\Box$

\vskip 0.5cm

\noindent
{\bf Acknowledgments}\\

\noindent The research of Deli Li was partially supported by a grant from the Natural Sciences and 
Engineering Research Council of Canada (grant \#: RGPIN-2014-05428)and the research of Han-Ying Liang 
was partially supported by the National Natural Science Foundation of China (grant \#: 11271286).

\vskip 0.5cm

{\bf References}

\begin{enumerate}

\item Adler, A., Rosalsky, A., Taylor, R. L.: Strong laws of large numbers for weighted sums
of random elements in normed linear spaces. Internat. J. Math. \& Math. Sci. {\bf 12}, 507-529
(1989).

\item de Acosta, A.: Inequalities for {\it B}-valued random vectors
with applications to the law of large numbers. Ann. Probab. {\bf 9},
157-161 (1981).

\item Feller, W.: A limit theoerm for random variables with infinite
moments. Amer. J. Math. {\bf 68}, 257-262 (1946).
        
\item Kolmogoroff, A.: Sur la loi forte des grands nombres. C. R.
Acad. Sci. Paris S\'{e}r. Math. {\bf 191} (1930), 910-912.

\item Ledoux, M., Talagrand, M.:  Probability in Banach Spaces:
Isoperimetry and Processes. Springer-Verlag, Berlin (1991).

\item Li, D., Liang, H.-Y., Rosalsky, A.: A note on symmetrization procedures 
for the laws of large numbers. Statist. Probab. Lett.  {\bf 121} (2017 a), 136-142.

\item Li, D., Liang, H.-Y., Rosalsky, A.: A probability inequality for sums of 
independent Banach space valued random variables. arXiv:1703.07868, (2017 b), 10 pages.

\item Marcinkiewicz, J., Zygmund, A.: Sur les fonctions
ind\'{e}pendantes. Fund. Math. {\bf 29} (1937), 60-90.

\item Martikainen, A. I., Petrov, V. V.: On a theorem of Feller.  Teor.
Veroyatnost. i Primenen. {\bf 25} (1980), ~194-197 (in Russian).
English translation in Theory Probab. Appl. {\bf 25} (1980), 191-193.

\end{enumerate}

\end{document}